\documentclass{amsart}

\usepackage{amsmath}
\usepackage{amssymb}
\usepackage[dvips]{color}

\begingroup 

\newtheorem{theorem}{Theorem}  
\newtheorem{lemma}[theorem]{Lemma}         
\endgroup

\newtheorem{definition}[theorem]{Definition}   

\newtheorem{example}[theorem]{Example}        

\numberwithin{equation}{section}

\newcommand{\sign}{\operatorname{sgn}}

\newcommand{\R}{{\mathbb R}}
\newcommand{\N}{{\mathbb N}}

\newcommand{\vare}{\varepsilon}

\newcommand{\lin}{\operatorname{span}}

\newcommand{\vp}{\varphi}

\newcommand{\diam}{\operatorname{diam}}
\newcommand{\card}{\operatorname{card}}

\newcommand{\mcI}{{\mathcal{I}}}

\newcommand{\mcM}{{\mathcal{M}}}

\newcommand{\mcP}{{\mathcal{P}}}

\begin{document}


\title{Second order differentiability of paths 
via a~generalized $\frac{1}{2}$-variation}

\author{Jakub Duda}
\email{duda@karlin.mff.cuni.cz}

\address{
Charles University, Department of Mathematical Analysis,
Sokolovsk\'a 83, 186 75 Praha 8 - Karl\'\i n, Czech Republic
}

\date{July 4, 2006}

\begin{abstract} 
We find an equivalent condition for a continuous vector-valued path  
to be Lebesgue equivalent to
a twice differentiable function.
For that purpose, we introduce the notion of a $VBG_{\frac{1}{2}}$ function,
which plays an analogous r\^ole
for the second order differentiability as the classical notion of a $VBG_*$ function
for the first order differentiability. 
In fact, for a function $f:[a,b]\to X$, being Lebesgue equivalent to
a twice differentiable function is the same as being Lebesgue equivalent to
a differentiable function with a pointwise Lipschitz derivative.
We also consider the case when the first derivative
can be taken non-zero almost everywhere.
\end{abstract}



\section{Introduction} 

Zahorski~\cite{Zah1} and Choquet~\cite{Ch} (see also Tolstov~\cite{Tol1})
proved a result
characterizing curves ($f:[a,b]\to\R^n$) that
allow a differentiable parametrization (resp.\ a dif.\ parametrization with almost everywhere non-zero
derivative) as those curves having
the $VBG_*$ property (resp.\ which are also not constant on any interval). 
Fleissner and Foran~\cite{FF} reproved this later (for real functions only and not considering
the case of a.e.\ nonzero derivatives) using a different result
of Tolstov.
The definition of $VBG_*$ is classical; see e.g.~\cite{S}.
The mentioned results were generalized by L.~Zaj\'\i\v{c}ek and the author~\cite{DZ}
to curves with values in Banach spaces
(and also metric spaces using the metric derivative instead
of the usual one).
Laczkovich, Preiss~\cite{LP}, and Lebedev~\cite{Leb}
studied (among other things) the case of $C^n$-parametrizations of
real-valued functions ($n\geq2$).
For a nice survey of differentiability of real-valued
functions via homeomorphisms, see~\cite{GNW}.
L.~Zaj\'\i\v{c}ek and the author~\cite{DZC2} characterized
the situation when a Banach space-valued curve  
admits a $C^2$-parametrization (for Banach spaces with a $C^1$ norm)
or a parametrization with finite convexity (for arbitrary Banach spaces).
\par
Let $X$ be a normed linear space, and $f:[a,b]\to X$. We say that $f$ is {\em Lebesgue equivalent} to $g:[a,b]\to X$
provided there exists a homeomorphism $h$ of $[a,b]$ onto itself such that $g=f\circ h$.
In the present note, we prove the following two theorems characterizing the
situation when a vector-valued path allows a~twice differentiable parametrization 
(resp.\ such a~parametrization with almost everywhere non-zero derivative):

\begin{theorem}\label{mainthm} 
Let $X$ be a normed linear space, and $f:[a,b]\to X$ be continuous.
Then the following are equivalent.
\begin{enumerate}
\item $f$ is Lebesgue equivalent to a twice differentiable function $g$.
\item $f$ is Lebesgue equivalent to a differentiable function $g$ whose
derivative is pointwise Lipschitz.
\item $f$ is $VBG_{\frac{1}{2}}$.
\end{enumerate}
\end{theorem}

\begin{theorem}\label{nonzthm} 
Let $X$ be a normed linear space, and $f:[a,b]\to X$ be continuous.
Then the following are equivalent.
\begin{enumerate}
\item $f$ is Lebesgue equivalent to a twice differentiable function $g$ with $g'(x)\neq0$ for a.e.\ $x\in[a,b]$.
\item $f$ is Lebesgue equivalent to a differentiable function $g$ whose
derivative is a pointwise Lipschitz function which is non-zero a.e.\ in $[a,b]$.
\item $f$ is $VBG_{\frac{1}{2}}$, 
and $f$ is not constant in any interval.
\end{enumerate}
\end{theorem}

As a matter of fact, a definition of a new notion of a $VBG_{\frac{1}{2}}$ function
(see Definition~\ref{vbgdef} below) involving 
a certain fractional variation, that was inspired by the results of 
Laczkovich, Preiss, and Lebedev, is necessary to achieve our goal.
\par
The case of $n$-times differentiable functions for $n\geq3$ is more
complicated even in the case $X=\R$, and this case
is treated in a separate paper~\cite{Dhigh}
(where we also prove a version of Zahorski lemma for $n$-times differentiable
homeomorphisms). 
The difficulty in the case of higher order derivatives of paths
stems from the fact that although for a curve parametrized by the arc-length,
the first derivative (provided it exists) is equal to the tangent (and
thus has norm~$1$), the magnitude of higher-order derivatives
is not thus simply bounded.
The proof in the real-valued case of $n\geq3$ uses some auxiliary variations and proceeds
in a rather indirect way.
This is a similar phenomenon as the case of $C^1$ parametrizations being
different from the case of $C^n$ ($n>1$) parametrizations; see e.g.~\cite[p.\ 405]{LP}
(since, in some sense, the $C^1$ case corresponds to twice-differentiable function
case).

\section{Preliminaries}

By $\lambda$ we will denote the Lebesgue measure on $\R$.
By $X$, we will always denote a normed linear space,
and by $B(x,r)$ an open ball with center $x$ and radius $r>0$.
If $X$ is separable, then it is well known that 
$X$ admits an equivalent G\^ateaux differentiable norm (see e.g.\ \cite{DGZ}).
For $f:[a,b]\to X$ we define the derivative $f'$ as usual (at the 
endpoints, we take the corresponding unilateral derivatives).
Similarly, the second derivative $f''(x)$ of $f$ at $x$ is defined
as $f''(x):=(f')'(x)$. Note that the property of ``being twice differentiable''
is preserved under equivalent renormings of $X$.
\par
We say that $f$ is {\em pointwise-Lipschitz at $x\in[a,b]$} provided
$\overline{\lim}_{t\to0} \frac{\|f(x+t)-f(x)\|}{|t|}$ is finite.
We say that $f$ is {\em pointwise-Lipschitz} provided
$f$ is pointwise-Lipschitz at each $x\in[a,b]$.
\par
Let $f:[a,b]\to X$ be continuous,
and assume that $X$ has a G\^ateaux differentiable norm
(there is no loss of generality
in this assumption since the continuity of
$f$ implies that $\overline{\lin}(f([a,b]))$ is separable).
By $K_f$ we will denote 
the set of points $x\in[a,b]$ such that there is no open interval~$U$ containing~$x$
such that $f|_{\overline{U}}$ is either constant or admits an arc-length parametrization
which is twice differentiable. 
\par
In the case of $X=\R$, the set $K_f$ coincides with the set of points of
varying monotonicity of $f$ (see e.g.\ \cite{LP}).
Obviously, $K_f$ is closed and $\{a,b\}\subset K_f$.
We easily see that $K_f$ does not depend on the choice of
the (equivalent) G\^ateaux smooth norm on $X$.
It is easy to see that if $f:[a,b]\to X$ is twice differentiable
(and $X$ has a G\^ateaux differentiable norm), then $f'(x)=0$ for each $x\in K_f$
by the chain rule for derivatives and by the continuity of $f'$.
\par
Let $K\subset[a,b]$ be a closed set with $a,b\in K$.
We say that an interval $(c,d)\subset[a,b]$ is {\em contiguous to $K$ in $[a,b]$}
provided $c,d\in K$ and $(c,d)\cap K=\emptyset$ (i.e.\ it is a maximal open
component of $[a,b]\setminus K$ in $[a,b]$).
\par
By $V(f,[a,x])$ we denote the (usual) variation of $f$ on $[a,x]$. We
will sometimes use the notation $v_f(x):=V(f,[a,x])$ for $x\in[a,b]$.
We say that $\{y_i\}^N_{i=0}$ is a {\em partition} of $[a,b]$ provided
$a=y_0<y_1<\dots<y_N=b$.
\par
We shall need the following  lemma. For a proof, see e.g.\ \cite[Lemma~2.7]{DZ}.

\begin{lemma}\label{basiclem} 
Let
$\{a,b\}\subset B\subset[a,b]$ be closed, and $f:[a,b]\to\R$ be
continuous. If $\lambda(f(B))=0$,
then we have
$V(f,[a,b])=\sum_{i\in{\mathcal I}}V(f,[c_i,d_i])$,
where $I_i=(c_i,d_i)$, $(i\in{\mathcal I}\subset\N)$
are all intervals contiguous to $B$ in $[a,b]$.
\end{lemma}

As in~\cite{LP}, for $g:[a,b]\to\R$, $\alpha\in(0,1)$, and $K\subset [a,b]$, we will define
$V_{\alpha}(g,K)$ as a supremum of sums
\[ \sum^m_{i=1} |g(b_i)-g(a_i)|^{\alpha},\]
where the supremum is taken over all collections
$\{[a_i,b_i]\}^m_{i=1}$ of non-overlapping intervals in $[a,b]$ with $a_i,b_i\in K$
for $i=1,\dots,m$.
\par
We will need the following auxiliary lemma:

\begin{lemma}\label{halftonullem}
Let $\alpha\in(0,1)$, $A\subset\R$ be bounded,
$f:A\to\R$ be uniformly continuous with
$V_{\alpha}(f,A)<\infty$. Then $\lambda(\overline{f(A)})=0$.
\end{lemma}

\begin{proof} 
By~\cite[Theorem~2.10]{LP} it follows that $SV_\alpha(f,A)=0$
(see~\cite{LP} for the definition of $SV_\alpha$).
It is easy to see that $SV_{\alpha}(f,A)=0$ implies $SV_1(f,A)=0$, and thus~\cite[Theorem~2.9]{LP}
shows that $\lambda(\overline{f(A)})=0$.
\end{proof}

We will need the following notion which plays the r\^{o}le of $VBG_*$ for
the second order differentiability.

\begin{definition}\label{vbgdef}
We say that a continuous $f:[a,b]\to X$ is {\em $VBG_{\frac{1}{2}}$} provided
$f$ has bounded variation, and
there exist closed sets $A_m\subset[a,b]$ $(m\in\mcM\subset\N)$
such that $K_f=\bigcup_{m\in\mcM} A_m$, 
and $V_{\frac{1}{2}}(v_f,A_m)<\infty$
for each $m\in\mcM$.
\end{definition}

It is easy to see that if $f$ is $VBG_{\frac{1}{2}}$ and $g$ is Lebesgue equivalent to $f$,
then $g$ is $VBG_{\frac{1}{2}}$. Also, it is easily seen that the class of $VBG_{\frac{1}{2}}$
functions does not depend on the equivalent norm of $X$.
\par
The following example shows that we cannot equivalently
replace $v_f$ by $f$ in Definition~\ref{vbgdef} (even in the case $X=\R$).

\begin{example} 
There exists a continuous function $f:[0,1]\to\R$ with bounded variation
such that $f$ is not $VBG_{\frac{1}{2}}$, but there exist closed $A_m\subset K_f$
such that $K_f=\bigcup_m A_m$, and
$V_{\frac{1}{2}}(f,A_m)<\infty$.
\end{example}

\begin{proof} 
Let $C\subset[0,1]$ be the standard middle-thirds Cantor set. 
By $\mcI_n$ we will denote the collection of all intervals
contiguous to $C$ such that $\lambda(I)<3^{-n}$ for $I\in\mcI_n$,
and by $K^n_{i}$, where $i=1,\dots,2^n$, $n\in\N$, denote
the closed intervals at level $n+1$ of the construction.
It is easy to see that there exist open intervals $I_{nik}\subset[0,1]$
and numbers $a_{ink}>0$, where $n,k\in\N$ and $i=1,\dots,2^n$,
such that
\begin{enumerate}
\item $I_{nik}\cap I_{n'i'k'}=\emptyset$ whenever $(n,i,k)\neq(n',i',k')$,
\item $\sum_{n,k\in\N}\sum_{i=1}^{2^n} a_{nik}<\infty$,
\item $\sum_{k\in\N}\sqrt{a_{nik}}=\infty$ whenever $n\in\N$ and $i=1,\dots,2^n$,
\item $\card\big\{(n,i,k):I_{nik}\subset I\big\}<\infty$ for all $m\in\N$ and $I\in\mcI_m$,
\item if $k\neq k'$, then there exists $x\in C$ such that either $I_{nik}<x<I_{nik'}$
or $I_{nik'}<x<I_{nik}$,
\item $I_{nik}\subset\big(K^n_i\cap\bigcup\mcI_n\big)$ for all $i=1,\dots,2^n$, $n,k\in\N$.
\end{enumerate}
Let $I=(a,b)\subset[0,1]$ be an open interval. We denote $l(I)=a$, $r(I)=b$, and $c(I)=\frac{a+b}{2}$.
We will define $f(x):=0$ when $x\in[0,1]\setminus\big(\bigcup_{n,k\in\N}\bigcup_{i=1}^{2^n}I_{nik}\big)$,
$f(c(I_{nik})):=a_{nik}$, and $f$ to be continuous and affine on $[l(I_{nik}),c(I_{nik})]$
and $[c(I_{nik}),r(I_{nik})]$. Then $f$ is a continuous function and by~(ii) it is easy
to see that $V(f,[0,1])<\infty$.
Index the countable family of closed sets 
\[\{C\}\cup\{\{l(I_{nik}),c(I_{nik}),r(I_{nik})\}:n,k\in\N,i=1,\dots,2^n\}\]
as $(A_m)_{m\in\N}$. It is easy to see that $K_f=\bigcup_{m\in\N}A_m$ and $V_{\frac{1}{2}}(f,A_m)<\infty$
for all $m\in\N$ (since $f|_C\equiv0$, and all those $A_m$ that satisfy $A_m\neq C$ are finite).
\par
Now we will show that $f$ is not $VBG_{\frac{1}{2}}$.
Suppose that $\tilde A_m$ satisfy $V_{\frac{1}{2}}(v_f,\tilde A_m)<\infty$,
and $K_f=\bigcup_m\tilde A_m$. Since $C=\bigcup_m(C\cap\tilde A_m)$,
by the Baire category theorem, there exists $m_0$ and an open
interval $U$ such that $C\cap U\subset C\cap\tilde A_{m_0}\cap U$ and $C\cap U\neq\emptyset$.
Thus, there exists $n\in\N$ and $i\in\{1,\dots,2^n\}$ such that
$K^n_i\subset U$, and conditions (v),(vi) imply that
\[ V_{\frac{1}{2}}(v_f,\tilde A_{m_0})\geq\sum_{\{I\in\mcI_n:I\subset K^n_i\}} 
\big({V(v_f,I)}\big)^{\frac{1}{2}}\geq\sum_k \sqrt{a_{nik}}=\infty,\]
which contradicts the choice of the sets $\tilde A_{m}$.
Thus, $f$ is not $VBG_{\frac{1}{2}}$.
\end{proof}

\section{Lemmata}

The following lemma is a sufficient condition for a function to be $VBG_{\frac{1}{2}}$.

\begin{lemma}\label{vbglem} 
Let $f:[a,b]\to X$ have a pointwise-Lipschitz derivative. Then $f$ is $VBG_{\frac{1}{2}}$.
\end{lemma}

\begin{proof}
Because $f'$ is continuous on $[a,b]$ (and thus bounded), we see that $f$
is Lipschitz (and thus has finite variation).
For $j\in\N$ define 
\[ D_j=\{x\in[a,b]:\|f'(x)-f'(z)\|\leq j|x-z|\text{ for all }z\in B(x,1/j)\}.\]
It is easy to see that $[a,b]=\bigcup_j D_j$, and $D_j$ is closed.
Let $D_j=\bigcup_{k\in\N} D_{jk}$ be such that each $D_{jk}$ is closed, and $\diam(D_{jk})<1/j$.
We order the doubly-indexed sequence $(K_f\cap D_{jk})_{j,k}$ into a single sequence
(while omitting empty sets); we will call the new sequence $A_m$ ($m\in\mcM\subset\N$). 
\par
It remains to show that $V_{\frac{1}{2}}(v_f,A_m)<\infty$, where $m\in\mcM$.
Let $m\in\mcM$, and fix $j,k\in\N$ such that $A_m=D_{jk}\cap K_f$.
Let $x<y$ be such that $x,y\in A_m$.
Note that
\begin{equation}\label{forind} 
|v_f(y)-v_f(x)|\leq\int^y_x \|f'(s)\|\,ds\leq j(y-x)^2.
\end{equation}
Applying~\eqref{forind} to $[x,y]=[a_i,b_i]$,
$i\in\{1,\dots,N\}$, where $[a_i,b_i]$ are non-overlapping intervals with $a_i,b_i\in A_m$,
we obtain
\begin{equation}\label{sumover}
\sum_{i=1}^N |v_f(b_i)-v_f(a_i)|^{\frac{1}{2}}\leq \sqrt{j}\sum_{i=1}^N
(b_i-a_i)\leq\sqrt{j}(b-a).
\end{equation}
By taking a supremum over all sequences $\{[a_i,b_i]\}^N_{i=1}$ as above,
we obtain that $V_{\frac{1}{2}}(v_f,A_m)<\infty$.
\end{proof}

\begin{lemma}\label{sqrlem}
Let $\zeta:[\sigma,\tau]\to\R$ be a continuous strictly increasing Lipschitz function with
$\zeta(0)=0$, and $\lambda(F)=0$ for some closed $F\subset[\sigma,\tau]$
with $\sigma,\tau\in F$. Then $\lambda(\sqrt{\zeta}(F))=0$, where $\sqrt{\zeta}(x):=\sqrt{\zeta(x)}$
for $x\in[\sigma,\tau]$.
\end{lemma}

\begin{proof}
Since the function $g(x)=\sqrt{x}$ on $[0,\infty)$ has property~(N)
(i.e.\ it maps zero sets onto zero sets), the conclusion easily follows.
\end{proof}

We will need the following simple lemma.

\begin{lemma}\label{homeolem}
Let $h_m:[a,b]\to[c_m,d_m]$ $(m\in\mcM\subset\N)$ be continuous increasing
functions such that $\sum_{m\in\mcM}h_m(x)<\infty$ for all $x\in[a,b]$.
Let $K\subset[a,b]$ be closed and such that $\lambda(h_m(K))=0$ for all $m\in\mcM$.
Then $h:[a,b]\to[c,d]$, defined as $h(x):=\sum_{m\in\mcM}h_m(x)$, 
is a continuous and increasing function $($for some
$c,d\in\R)$ such that $\lambda(h(K))=0$.
\end{lemma}

\begin{proof}
The continuity and monotonicity of $h$ follows easily by the assumptions.
Let $K\subset[a,b]$ be closed with $\lambda(h_m(K))=0$ for all $m\in\mcM$.
Without any loss of generality, we can assume that $\{a,b\}\subset K$.
Let $(c_p,d_p)$ ($p\in\mcP\subset\N$) be all the intervals contiguous to $K$ in~$[a,b]$.
Let $\vare>0$ and find $M\in\N$ such that $\sum_{m\in\mcM\cap(M,\infty)}(h_m(b)-h_m(a))<\vare$.
Then
\begin{equation*}
\begin{split}
\lambda(h([a,b]))&=\sum_{m\in\mcM}(h_m(b)-h_m(a))
\leq\vare+\sum_{m\in\mcM\cap[1,M]}(h_m(b)-h_m(a))\\
&=\vare+\sum_{m\in\mcM\cap[1,M]}\sum_{p\in\mcP}(h_m(d_p)-h_m(c_p))
\leq \vare+\sum_{p\in\mcP}\lambda(h(c_p,d_p)),
\end{split}
\end{equation*}
where we used Lemma~\ref{basiclem} to obtain the second equality. 
Since $\card\big( h((c_p,d_p))\cap h((c_q,d_q))\big)\leq 1$ for $p,q\in\mcP$,
$p\neq q$,
we obtain the equality 
\[\lambda(h([a,b]))=\lambda\big(h\big(\bigcup_{p\in\mcP}(c_p,d_p)\big)\big).\]
Since the set $h(K)\cap h\big(\bigcup_{p\in\mcP}(c_p,d_p)\big)$ is countable,
we get $\lambda(h(K))=0$.
\end{proof}

\begin{lemma}\label{varlem}
Suppose that $X$ is a normed linear space with
a G\^ateaux smooth norm.
Let $f:[a,b]\to X$ be a continuous $VBG_{\frac{1}{2}}$ function
which is not constant on any interval.
Then there exists a continuous strictly increasing $v:[a,b]\to[\alpha,\beta]$ 
such that $\lambda(v(K_f))=0$, $f\circ v^{-1}$ is twice differentiable
on~$[\alpha,\beta]\setminus v(K_f)$ with $(f\circ v^{-1})'(x)\neq0$
for $x\in[\alpha,\beta]\setminus v(K_f)$, and
for each $x\in K_f$ there exists $0<C_x<\infty$ such that
\begin{equation}\label{vartimes} 
\|f(y)-f(z)\|\leq C_x |v(z)-v(y)|(|v(z)-v(x)|+|v(y)-v(x)|),
\end{equation}
whenever $y,z\in[a,b]$, and $\sign(y-x)=\sign(z-x)$.
\end{lemma}

\begin{proof} 
Let $A_m$ ($m\in\mcM\subset\N$) be as in the definition of $VBG_{\frac{1}{2}}$ for $g=f\circ v_f^{-1}$.
Note that $g$ is $1$-Lipschitz, and $K_g=v_f(K_f)$.
Since $f$ is $VBG_{\frac{1}{2}}$, by Lemma~\ref{halftonullem} we have
$\lambda(v_f(K_f))=\lambda(v_g(K_g))=0$.
Let $\ell=v_f(b)$. 
Note that because $g$ is an arc-length parametrization
of $f$, we have $V(g,[c,d])=d-c$ for all $0\leq c<d\leq\ell$
(we will use this fact frequently without necessarily repeating it).
Let $(c_p,d_p)$ ($p\in\mcP\subset\N$) be all the intervals contiguous to $K_g$ in $[0,\ell]$.
Since $\lambda(v_g(K_g))=0$, by Lemma~\ref{basiclem} (applied to $f=v_g$) we have
$V(g,[0,\ell])=\ell=\sum_{p\in\mcP}V(g,[c_p,d_p])=\sum_{p\in\mcP}(d_p-c_p)$,
and thus $\lambda(K_g)=\ell-\lambda\big(\bigcup_{p\in\mcP}(c_p,d_p)\big)=0$.
For $m\in\mcM$ and $x\in[0,\ell]$, we define $v_m(x)$ as a supremum of the sums
\begin{equation}\label{vmdef}
\sum_{i=1}^{N} ( b_i-a_i)^{\frac{1}{2}},
\end{equation}
where the supremum is taken over all finite sequences $\{[a_i,b_i]\}^N_{i=1}$
of non-overlapping intervals in $[0,\ell]$ such that 
$a_i,b_i\in (A_m\cup\{0,x\})\cap[0,x]$ for $i=1,\dots,N$.
Similarly, we define $\tilde v_m(x)$ for $x\in[0,\ell]$ as a supremum of the sums in~\eqref{vmdef},
where the supremum is taken over all finite sequences 
$\{[a_i,b_i]\}^N_{i=1}$
of non-overlapping intervals in $[0,\ell]$ such that 
$a_i,b_i\in (A_m\cup\{x,\ell\})\cap[x,\ell]$ for $i=1,\dots,N$.
Note that $v_m$ is increasing and $\tilde v_m$ is decreasing on $[0,\ell]$.
Note that $g$ is affine on each $[c_p,d_p]$, and 
\begin{equation}\label{vmequa}
v_m(x)=v_m(z)+(x-z)^{\frac{1}{2}}\quad\text{ for }x\in[c_p,d_p],
\end{equation}
where $z=\max((A_m\cup\{0\})\cap[0,c_p])$, and similarly for $\tilde v_m$.
Thus $v_m$ (and similarly~$\tilde v_m$) is twice (or even infinitely many times) 
differentiable on~$[0,\ell]\setminus v_f(K_f)$ with $v_m'(x)>0$ for all
$x\in[0,\ell]\setminus v_f(K_f)$.
Find $\vare_m>0$ such that 
\begin{itemize}
\item[(a)] if we define $w(x):=\sum_m \vare_m\cdot (v_m(x)-\tilde v_m(x))$,
then $w(0)$, and $w(\ell)$ are finite (and thus $w(x)$ is finite for all $x\in[0,\ell]$),
and $w$ is continuous on $[0,\ell]$ (provided all $v_m,\tilde v_m$ were continuous),
\item[(b)] for all $m\in\mcM$ and $p\in\mcP$ with $c_p+1/m<d_p-1/m$ and all $x\in(c_p+1/m,d_p-1/m)$,
we have $\vare_m\cdot\max(v_m'(x),|v_m''(x)|,-\tilde v'_m(x),|\tilde v''_m(x)|)<2^{-m}$.
\end{itemize}
By~(b), it is easy to see that $w'(x)$ exists, is positive, and $w''(x)$ exists for each $x\in[0,\ell]\setminus v_f(K_f)$.
Put $v:=w\circ v_f$, $\alpha=v(a)$, and $\beta=v(b)$.
\par
To show that $v$ is strictly increasing, it is enough to show that
$w$ is strictly increasing (as $v_f$ is strictly increasing by the
fact that $f$ is not constant on any interval).
On the other hand, to show that $w$ is strictly increasing, it is enough
to show that $v_m$ is strictly increasing for each $m\in\mcM$.
Fix $m\in\mcM$. 
Let $x,y\in[0,\ell]$ with $x<y$. If $x,y\in[c_p,d_p]$ for
some $p\in\mcP$, then~\eqref{vmequa} implies that $v_m(x)<v_m(y)$,
and similarly if $x\in(c_p,d_p)$ or $y\in(c_{p'},d_{p'})$ for some 
$p\in\mcP$ (resp.\ $p'\in\mcP$). If $x,y\in K_f$, and $(x,y)\cap A_m=\emptyset$,
then 
\begin{equation}\label{vmeq1}
v_m(t)=v_m(z)+\sqrt{t-z}\quad\text{ for all }t\in[x,y],
\end{equation}
where $z=\max((A_m\cup\{0\})\cap[0,x])$, and thus $v_m(x)<v_m(y)$.
Finally, if there exists $q\in A_m\cap (x,y)$, then
$v_m(x)\leq v_m(q)<v_m(q)+\sqrt{y-q}\leq v_m(y)$,
and thus $v_m(x)<v_m(y)$ also in this case.
By a similar argument, $\tilde v_m$ is strictly decreasing.
\par
For a fixed $m\in\mcM$, we will prove that whenever $r,s\in A_m\cup\{0,\ell\}$
with $r<s$, then
\begin{equation}\label{vmprop} 
v_m(s)-v_m(r)\leq \sum_{\substack{p\in\mcP:\\(c_p,d_p)\cap[r,s]\neq\emptyset}}
(v_m(d_p)-v_m(c_p)).
\end{equation}
A symmetrical argument then shows that
\begin{equation}\label{vmproptilde} 
\tilde v_m(r)-\tilde v_m(s)\leq \sum_{\substack{p\in\mcP:\\(c_p,d_p)\cap[r,s]\neq\emptyset}}
(\tilde v_m(c_p)-\tilde v_m(d_p)).
\end{equation}
To prove~\eqref{vmprop}, fix $\vare_0>0$, and
let $\{[a_i,b_i]\}^{N}_{i=1}$ be non-overlapping intervals
in $[r,s]$ such that $a_i,b_i\in (A_m\cup\{r,s\})\cap[r,s]$ for $i=1,\dots,N$
such that
$v_m(s)=v_m(r)+\sum_{i=0}^{N-1} (b_{i}-a_i)^{\frac{1}{2}}+\vare$,
for some $0\leq\vare<\vare_0/2$.
For $i\in\{1,\dots,N\}$ by Lemma~\ref{basiclem} applied to
$f=g$ on $[a,b]=[a_i,b_{i}]$ and $B=(A_m\cup\{r,s\})\cap[a_i,b_{i}]$ (note that
$\lambda(g(A_m))=0$ since $\lambda(g(K_g))=0$, and thus $\lambda(g(B))=0$), 
let $(\gamma^i_j,\delta^i_j)$ ($j\in\{1,\dots,J^i\}$) 
be a finite collection of intervals contiguous to $A_m\cup\{r,s\}$ in $[a_i,b_{i}]$ such that
$(b_{i}-a_i)\leq \sum_{j=1}^{J^i} (\delta^i_j-\gamma^i_j)+\big(\frac{\vare_0}{2N}\big)^2$.
Then 
\begin{equation}\label{pee1} 
v_m(s)-v_m(r) \leq \sum_{i=1}^{N} \sum_{j=1}^{J^i} (\delta^i_j-\gamma^i_j)^{\frac{1}{2}}+\frac{\vare_0}{2}+\vare.
\end{equation}
By Lemma~\ref{sqrlem} applied to $\zeta(x)=x-\gamma^i_j$ on $[\sigma,\tau]=[\gamma^i_j,\delta^i_j]$,
$F=K_g\cap [\gamma^i_j,\delta^i_j]$, and
because $v_m(x)=v_m(\gamma^i_j)+(x-\gamma^i_j)^{\frac{1}{2}}$ for $x\in[\gamma^i_j,\delta^i_j]$,
we have that $\lambda(v_m(K_g\cap[\gamma^i_j,\delta^i_j]))=0$, and by Lemma~\ref{basiclem}
applied to $f=v_m$ on $[a,b]=[\gamma^i_j,\delta^i_j]$, and $B=K_g\cap[\gamma^i_j,\delta^i_j]$,
we obtain that 
$(\delta^i_j-\gamma^i_j)^{\frac{1}{2}} \leq \sum_{\substack{p\in\mcP:\\(c_p,d_p)\subset [\gamma^i_j,\delta^i_j]}} (v_m(d_p)-v_m(c_p))$
for each $i\in\{1,\dots, N\}$ and $j\in\{1,\dots,J^i\}$. 
Combining this inequality with~\eqref{pee1},
we get
$ v_m(s)-v_m(r)\leq \sum_{\substack{p\in\mcP:\\(c_p,d_p)\cap[r,s]\neq\emptyset}}
(v_m(d_p)-v_m(c_p))+\vare_0$,
and by sending $\vare_0\to0$ it follows that~\eqref{vmprop} holds.
\par
To show that $v$ is continuous, it is enough to show that each $v_m$
is continuous (as this implies that $w$ is continuous by the choice of $\vare_m$'s, and the continuity of~$v_f$ follows from e.g.\ \cite[\S 2.5.16]{Fe}).
Fix $m\in\mcM$. From~\eqref{vmequa}, it follows that 
\begin{itemize}
\item[$(*)$]
$v_m$~is continuous from the right at all points $x\in\bigcup_{p\in\mcP}[c_p,d_p)$,
and continuous from the left at all points $x\in\bigcup_{p\in\mcP}(c_p,d_p]$.
\end{itemize}
If $(x,y)\cap A_m=\emptyset$ for some $y>x$ with $y\in(0,\ell]\cap K_g$, 
then~\eqref{vmeq1} implies that $v_m$ is continuous from the right at $x$.
If $x\in A_m$ is a right-hand-side accumulation point of $A_m$ (i.e.\ $A_m\cap(x,x+\delta)\neq\emptyset$
for all $\delta>0$), then~\eqref{vmprop} implies that 
$\lim_{\substack{y\to x+\\y\in A_m}}v_m(y)=v_m(x)$,
since
\begin{equation}\label{tozero}
\sum_{\substack{p\in\mcP\\(c_p,d_p)\cap[x,y]\neq\emptyset}}(v_m(d_p)-v_m(c_p))\to0
\end{equation}
as $y\to x+$. Now the monotonicity of $v_m$ implies that it is continuous from
the right at~$x$. Concerning the continuity from the left, by~$(*)$ it is enough
to prove that $v_m$ is continuous from the left at all points $y\in(K_g\cap(0,\ell])\setminus
\bigcup_{p\in\mcP}\{d_p\}$. Fix such a point~$y$. If there is an $x\in[0,y)$ such that 
$(x,y)\cap A_m=\emptyset$, then~\eqref{vmeq1} implies that $v_m$ is continuous from the
left at~$y$. If~$y$ is a left-hand-side
accumulation points of~$A_m$, then~\eqref{vmprop} together
with~\eqref{tozero} imply that~$v_m$ is continuous from the left at~$y$.
A similar argument as above yields the continuity of $\tilde v_m$.
\par
Now we will prove that $\lambda(v(K_f))=0$. Note that we already established
that $\lambda(K_g)=0$.
Because $K_g=v_f(K_f)$, it is enough to prove that $\lambda(w(K_g))=0$.
To apply Lemma~\ref{homeolem} to $h_k$, where $h_{2k}:=\vare_k\cdot v_k$, and $h_{2k+1}:=-\vare_k\cdot\tilde v_k$, we have to check that $\lambda(v_m(K_g))=0$ and $\lambda(\tilde v_m(K_g))=0$
for all $m\in\mcM$.
Let $m\in\mcM$. Then~\eqref{vmprop} applied to $r=0$, and $s=\ell$ shows that
$v_m(\ell)-v_m(0)\leq \sum_{p\in\mcP} (v_m(d_p)-v_m(c_p))$,
and since $v_m(K_g)\cap v_m\big(\bigcup_{p\in\mcP}(c_p,d_p)\big)=\emptyset$, 
we get $\lambda(v_m(K_g))=0$. Similarly, we obtain $\lambda(\tilde v_m(K_g))=0$.
Thus, Lemma~\ref{homeolem} shows that $\lambda(w(K_g))=0$.
\par
To prove that the second derivative of $f\circ v^{-1}$ exists
and the first derivative is non-zero on $[\alpha,\beta]\setminus v(K_f)$, let
$x\in[\alpha,\beta]\setminus v(K_f)$. Put $y=w^{-1}(x)$. 
There exists $p\in\mcP$ and $q\in\N$ such that $y\in(c_p+1/q,d_p-1/q)$.
Since (by the chain rule 
and the smoothness of the norm
on $X$) $g$ is twice differentiable on $(c_p,d_p)$ and
$\|g'(x)\|=1$ for all $x\in(c_p,d_p)$ (because $g$ is the arc-length
parametrization of $f$ and $g'$ is continuous on $(c_p,d_p)$), 
it is enough to prove that $w'(y)$ exists, is non-zero,
and $w''(y)$ exists (since then $(f\circ v^{-1})'(x)=g'(y)\cdot(w^{-1})'(x)$, and 
$(f\circ v^{-1})''(x)=g''(y)\cdot ((w^{-1})'(x))^2+g'(y)\cdot(w^{-1})''(x)$).
But by the choice of $\vare_m$ (for $m>q$), and by the properties of $v_m$, $\tilde v_m$ for all $m$, 
it is easy to see that $w'(y)$ exists, $w'(y)>0$, and
$w''(y)$ exists; the rest is a straightforward application of the ``derivative of
the inverse'' rule.
\par
To prove~\eqref{vartimes} for $f$ and $v$, by a substitution using $v_f$, it is easy
to see that it is enough to establish a version of~\eqref{vartimes}, where~$f$ is replaced by~$g$, and
$v$ by~$w$. 
To~that end, take $m\in\mcM$ such that $x\in A_m$, and let $C_m=(\vare_m)^{-2}$. 
Take $y,z\in[0,\ell]$. Without any loss of generality,
we can assume that $x<y<z$ (if $y<x$, then a symmetric estimate using $\tilde v_m$ yields the conclusion). 
Let $0<\vare_0<v_m(z)-v_m(x)$. 
Find a sequence $\{[a_i,b_i]\}_{i=1}^N$ of non-overlapping intervals with endpoints in $(A_m\cup\{x,y\})\cap[x,y]$ with $b_i<a_{i+1}$ for $i=1,\dots,N-1$, and 
such that $v_m(y)=v_m(x)+\sum_{i=1}^{N} (b_i-a_i)^{\frac{1}{2}}+\vare$,
for some $0\leq\vare\leq\vare_0$ (such a choice of $\{[a_i,b_i]\}^N_{i=1}$ is possible because $x\in A_m$).
Then
\begin{equation*}
\begin{split} 
v_m(z)&\geq  v_m(x)+\sum_{i=1}^{N-1} (b_i-a_i)^{\frac{1}{2}}+(z-a_N)^{\frac{1}{2}}\\
&\geq v_m(x)+\sum_{i=1}^{N-1} (b_i-a_i)^{\frac{1}{2}}+(z-y+b_N-a_N)^{\frac{1}{2}}
\end{split}
\end{equation*}
Thus (since $z-y\geq\|g(z)-g(y)\|$), we get
\begin{equation*}
\begin{split} 
v_m(z)-v_m(y)
&\geq  (b_N-a_N+z-y)^{\frac{1}{2}}-({b_N-a_N})^{\frac{1}{2}}-\vare\\
&=\frac{\|g(y)-g(z)\|}{(b_N-a_N+z-y)^{\frac{1}{2}}+(b_N-a_N)^{\frac{1}{2}}}-\vare\\  
&\geq\frac{ \|g(y)-g(z)\|}{v_m(z)-v_m(x)+v_m(y)-v_m(x)}-\vare.
\end{split}
\end{equation*}
By sending $\vare_0\to0$, we obtain
\begin{equation}\label{intm}
\|g(y)-g(z)\|\leq (v_m(z)-v_m(y))(v_m(z)-v_m(x)+v_m(y)-v_m(x)).
\end{equation} 
To finish the proof of~\eqref{vartimes} for $g$ and $w$,
note that $v_m(\tau)-v_m(\sigma)\leq \frac{1}{\vare_m} (w(\tau)-w(\sigma))$
for any $0\leq\sigma<\tau\leq\ell$;
thus~\eqref{vartimes} follows from~\eqref{intm}.
\end{proof}

We will need the following version of Zahorski's lemma. See e.g.~\cite{GNW}
for a proof of a slightly weaker statement.

\begin{lemma}\label{Zahorski}
Let $F\subset[\alpha,\beta]$ be closed, $\{\alpha,\beta\}\subset F$, and $\lambda(F)=0$. Then
there exists an $($increasing$)$ continuously differentiable homeomorphism $h$ of $[\alpha,\beta]$ 
onto itself such that $h'(x)=0$ if and only if $x\in h^{-1}(F)$, $h$~is
twice differentiable on $[\alpha,\beta]\setminus h^{-1}(F)$, and $h^{-1}$ is absolutely continuous.
\end{lemma}

\begin{proof} 
Since we were not able to locate a reference in the 
literature for this exact statement, we will
sketch the proof. 
Let $(a_i,b_i)$ (where $i\in\mcI\subset\N$) be all the intervals contiguous to $F$ in $[\alpha,\beta]$.
For each $i\in\mcI$ find a $C^1$ function $\psi_i:(a_i,b_i)\to\R$ such that
\begin{itemize}
\item $\psi_i(x)\geq0$ for all $x\in(a_i,b_i)$, and $\lim_{x\to a_i+}\psi_i(x)=\lim_{x\to b_i-}\psi_i(x)=\infty$,
\item $m_i:= \min_{x\in(a_i,b_i)}\psi_i(x)>0$, and if $|\mcI|=\aleph_0$,
then $\lim_{\substack{i\to\infty\\i\in\mcI}}m_i=\infty$, 
\item $\sum_{i\in\mcI}\int^{b_i}_{a_i}\psi_i(t)\,dt<\infty$.
\end{itemize}
Such functions $\psi_i$ clearly exist.
Define $\psi:[\alpha,\beta]\to\R$ as $\psi(x):=\psi_i(x)$ for $x\in(a_i,b_i)$, and $\psi(x)=0$ for $x\in F$.
It is easy to see that $\psi$ is integrable. 
Define $k(x):=\int^x_\alpha \psi(t)\,dt$; then $k$ is continuous and (strictly) increasing. 
By integrability of $\psi$, it follows that $k$ has Luzin's property~$(N)$, 
and thus $k$ is absolutely continuous by the Banach-Zarecki theorem (see e.g.~\cite[Theorem~3]{Var}).
It is easy to~see that $k$ is twice differentiable on $[\alpha,\beta]\setminus F$ with $k'(x)>0$.
We also have that $k'(x)=\infty$ for $x\in F\setminus\big(\bigcup_i\{a_i\}\big)$, as for $x\in F$ and $t>0$
small enough, we 
have
\[ k(x+t)-k(x)\geq m_j(x+t-a_j)+\sum_{(a_i,b_i)\subset[x,x+t]} m_i(b_i-a_i)\geq m_t\cdot t,\]
where $j\in\mcI$ is such that $x+t\in(a_j,b_j)$ and 
for $m_t:=\min\{m_k:(a_k,b_k)\cap[x,x+t]\neq\emptyset\}$ we have $\lim_{t\to0+}m_t=\infty$
by the choice of $\psi_i$. 
If $x=a_i$ for some $i\in\mcI$, then we have $k(x+t)-k(x)\geq t\cdot\min_{y\in[x,x+t]}\psi_i(y)$,
and the minimum goes to infinity with $t\to0+$ by the choice of~$\psi_i$.
By continuity and symmetry, the rest follows.
Now define $\vp(x):=\alpha+\frac{\beta-\alpha}{k(\beta)}k(x)$, $h:=\vp^{-1}$, and the lemma easily follows.

\end{proof}

\section{Proofs of the main results}

\begin{proof}[Proof of Theorem~\ref{mainthm}.]
The implication (i)$\implies$(ii) is trivial.
To prove that (ii)$\implies$(iii), let $h$ be a homeomorphism such that $g=f\circ h$ has pointwise-Lipschitz
derivative.
Then Lemma~\ref{vbglem} implies that $g$ is $VBG_{\frac{1}{2}}$.
By a~remark following Definition~\ref{vbgdef}, it follows that $f$ is $VBG_{\frac{1}{2}}$.
\par
To prove that (iii)$\implies$(i), 
without any loss of generality, we can assume that the norm on $X$
is G\^ateaux differentiable (since $\overline{\lin}(f([a,b]))$ is separable
and second order differentiability of a path does not depend on the equivalent
norm on $X$).
First, assume that $f$ is not
constant on any interval. Lemma~\ref{varlem} implies that there
exists an increasing homeomorphism $v:[a,b]\to[\alpha,\beta]$ such that $f\circ v^{-1}$ is 
differentiable on $[\alpha,\beta]$, twice
differentiable on $[\alpha,\beta]\setminus v(K_f)$, and 
$\lambda(v(K_f))=0$. Apply Lemma~\ref{Zahorski} to $F=v(K_f)$ 
to obtain an (increasing) continuously differentiable
homeomorphism $h:[\alpha,\beta]\to[\alpha,\beta]$ such that
$h'(x)=0$ iff $x\in h^{-1}(v(K_f))$, and such that $h$ is
twice differentiable on $[\alpha,\beta]\setminus h^{-1}(v(K_f))$.
Let $g=f\circ v^{-1}\circ h$. By the chain rule for derivatives,
we have that $g$ is twice differentiable on $[\alpha,\beta]\setminus h^{-1}(v(K_f))$.
Let $x\in h^{-1}(v(K_f))$. Then by~\eqref{vartimes} there exists a $C_x>0$ such that
\begin{equation}\label{vv1}
\frac{\|f\circ v^{-1}(y)-f\circ v^{-1}(z)\|}{|y-z|}\leq 2 C_x\,|z-h(x)|
\end{equation}
for $z<y<h(x)$ or $h(x)<y<z$ (and by continuity this holds also for $y=h(x)$), and $y,z\in[\alpha,\beta]$.
It follows that $(f\circ v^{-1})'(h(x))=0$. Thus $g'(x)=0$ by the chain rule.
It also follows from~\eqref{vv1} that $(f\circ v^{-1})'(\cdot)$ is pointwise-Lipschitz at $h(x)$ with constant $2C_x$.
This implies that
\begin{equation*}
\begin{split} 
\bigg\|\frac{g'(x+t)-g'(x)}{t}\bigg\|&=\bigg\|\frac{(f\circ v^{-1})'(h(x+t))h'(x+t)}{t}\bigg\|\\
&=\bigg\|\frac{(f\circ v^{-1})'(h(x+t))-(f\circ v^{-1})'(h(x))}{t}\bigg\|\cdot h'(x+t)\\
&\leq 2C_x\cdot\bigg|\frac{h(x+t)-h(x)}{t}\bigg|\cdot h'(x+t),
\end{split}
\end{equation*}
for all $x+t\in[\alpha,\beta]$. The continuity of~$h'$ at~$x$
shows that $g''(x)=0$. It is easy to~see that $f$ is Lebesgue equivalent to~$g$ (by composing $v^{-1}\circ h$ with an affine change of~parameter).
\par
If $f$ is constant on some interval, then let $(c_i,d_i)$ ($i\in\mcI\subset\N$) be the collection of all maximal open intervals such that $f$ is constant on each $(c_i,d_i)$.
It is easy to~see that we can find a continuous function $\tilde f:[a,b]\to X$ such that
$f=\tilde f$ on $[a,b]\setminus\bigcup_i(c_i,d_i)$, $\tilde f$ is affine and non-constant
on $(c_i,(c_i+d_i)/2)$, $((c_i+d_i)/2,d_i)$, and such that $\tilde f$ is $VBG_{\frac{1}{2}}$.
By the previous paragraph, there exists a homeomorphism~$h$ of~$[a,b]$ onto itself
such that $\tilde f\circ h$ is twice differentiable. It follows that $f\circ h$
is twice differentiable (since $(\tilde f\circ h)'(x)=(\tilde f\circ h)''(x)=0$ for all $x\in\bigcup_i\{c_i,d_i\}$
by the construction).
\end{proof}

\begin{proof}[Proof of Theorem~\ref{nonzthm}.] The implication (i)$\implies$(ii) is trivial.
To prove that (ii)$\implies$(iii), note that if $g'(x)\neq0$ for a.e.\ $x\in[a,b]$,
then $g$ is not constant in any interval. This notion is clearly stable with respect
to Lebesgue equivalence. The rest follows from Theorem~\ref{mainthm}.
\par
To prove that~(iii)$\implies$(i), we can follow the proof of the corresponding implication of Theorem~\ref{mainthm}.
To see that the resulting function $g$ has non-zero derivative almost everywhere, we note that
the homeomorphism~$h$ obtained by applying the Lemma~\ref{Zahorski} has an absolutely continuous inverse.
The rest follows easily.
\end{proof}

The following example shows that even in the case of $X=\R$,
$VBG_{1/2}$ functions do not coincide 
with continuous functions satisfying $V_{1/2}(f,K_f)<\infty$.

\begin{example}
There exists a continuous $VBG_{1/2}$ function $f:[0,1]\to\R$ such
that $V_{1/2}(f,K_f)=\infty$ $($and thus $f$ is not Lebesgue equivalent
to a $C^2$ function by~\cite[Remark~3.6]{LP}$)$.
\end{example}

\begin{proof}
Let $a_n\in(0,1)$ be such that $a_n\downarrow0$. Define
$f(a_{2k})=0$, $f(a_{2k+1})=1/k^2$ for $k=1,\dots$, and $f(0)=f(1)=0$.
Extend $f$ to be continuous and affine on the intervals $[a_{2k+1},a_{2k}]$
and $[a_{2k+2},a_{2k+1}]$. Then $K_f=\{0,1\}\cup\{a_n:n\geq2\}$ and 
it is easy to see that $f$ is $VBG_{1/2}$ but $V_{1/2}(f,K_f)=\infty$.
\end{proof}

\section*{Acknowledgment} The author would like to thank Professor Lud\v ek Zaj\'\i\v cek
for a valuable conversation that led to several simplifications and improvements of the paper.

\end{document}